\documentclass[12pt]{article}

\usepackage{amsfonts}

\usepackage{amsmath,amssymb}

 \numberwithin{equation}{subsection}

\begin{document}

\title{Properties of hitting times for $G$-martingale  }
\author{Yongsheng Song\\
\small Academy of Mathematics and Systems Science, \\
\small Chinese Academy of Sciences, Beijing, China;\\
\small yssong@amss.ac.cn}

\date{}

\maketitle

\begin{abstract}

In this article, we consider the properties of hitting times for
$G$-martingale and the stopped processes. We prove that the stopped
processes for $G$-martingales are still $G$-martingales and that the
hitting times for  a class of $G$-martingales including $G$-Brownian
motion are quasi-continuous. As an application, we improve the
$G$-martingale representation theorems in [Song10].

\end{abstract}

\maketitle
\section{Introduction }
Recently, [P06], [P08] introduced the notion of sublinear
expectation space, which is a generalization of probability space.
One of the most important sublinear expectation space is
$G$-expectation space. As the counterpart of Wiener space in the
linear case, the notions of $G$-Brownian motion, $G$-martingale, and
It$\hat{o}$ integral w.r.t $G$-Brownian motion were also introduced.
These notions have very rich and interesting new structures which
nontrivially generalize the classical ones.

 As is well known, stopping times play a great role in classical stochastic
analysis. However, it is difficult to apply stopping time technique
in subliner expectation space since the stopped process maynot
belong to the class of processes which are meaningful in the present
situation. For example, let $\{M_t\}_{t\in[0, T]}$ be a
$G$-martingale and $\tau$ be an $\mathbb{F}$-stopping time, we don't
know whether $M^\tau_t$ has a quasi-continuous version for $t\in[0,
T]$.

In this article  we consider the properties of hitting times for
$G$-martingale and the stopped processes. We prove that the stopped
processes for $G$-martingales are still $G$-martingales and that the
hitting times for  symmetric $G$-martingales with strictly
increasing quadratic variation processes are quasi-continuous. As an
application, we prove that  any symmetric random variable can be
approximated by bounded random variables that are also symmetric.
Besides, we improve the results in [Song10] for $G$-martingale
representation by a stopping time technique.

This article is organized as follows: In section 2, we recall some
basic notions and results of $G$-expectation and the related space
of random variables. In section 3, we give several preliminary
lemmas. In section 4, we prove that the stopped processes for
$G$-martingales are still $G$-martingales and that the hitting times
for  a class of $G$-martingales including $G$-Brownian motion are
quasi-continuous. In section 5, we give some applications by a
stopping time technique.

\section{Preliminary }
We recall some basic notions and results of $G$-expectation and the
related space of random variables. More details of this section can
be found in [P07].

\subsection{G-expectation }

\noindent {\bf Definition 2.1} Let $\Omega$
 be a given set and let ${\cal H}$ be a linear space of real valued
functions defined on $\Omega$
 with $c \in {\cal H}$ for all constants $c$. ${\cal H}$ is considered as the
space of  ¡°random variables¡±. A sublinear expectation $\hat{E}$ on
${\cal H}$ is a functional $\hat{E}: {\cal H}\rightarrow R $
satisfying the following properties: for all $X, Y \in {\cal H}$, we
have

(a) Monotonicity: If $X\geq Y$ then $\hat{E}(X) \geq \hat{E} (Y)$.

(b) Constant preserving: $\hat{E} (c) = c$.

(c) Sub-additivity: $\hat{E}(X)-\hat{E}(Y) \leq \hat{E}(X-Y)$.

(d) Positive homogeneity: $\hat{E} (\lambda X) = \lambda \hat{E}
(X)$, $\lambda \geq 0$.

\noindent$(\Omega, {\cal H}, \hat{E})$ is called a sublinear
expectation space.

\noindent {\bf Definition 2.2} Let $X_1$ and $X_2$ be two
$n$-dimensional random vectors defined respectively in sublinear
expectation spaces $(\Omega_1, {\cal H}_1, \hat{E}_1)$ and
$(\Omega_2, {\cal H}_2, \hat{E}_2)$. They are called identically
distributed, denoted by $X_1 \sim X_2$, if $\hat{E}_1[\varphi(X_1)]
= \hat{E}_2[\varphi(X_2)]$, $\forall \varphi\in C_{l, Lip}(R^n)$,
where $ C_{l, Lip}(R^n)$ is the space of real continuous functions
defined on $R^n$ such that $$|\varphi(x) - \varphi(y)| \leq C(1 +
|x|^k + |y|^k)|x - y|, \forall x, y \in R^n,$$ where $k$ depends
only on $\varphi$.

\noindent {\bf Definition 2.3} In a sublinear expectation space
$(\Omega, {\cal H}, \hat{E})$ a random vector $Y = (Y_1,
\cdot\cdot\cdot, Y_n)$, $Y_i \in {\cal H}$ is said to be independent
to another random vector $X = (X_1, \cdot\cdot\cdot, X_m)$, $X_i \in
{\cal H}$ under $\hat{E}(\cdot)$ if for each test function
$\varphi\in C_{l, Lip}(R^m\times R^n)$ we have $\hat{E}[\varphi(X, Y
)] = \hat{E}[\hat{E} [\varphi(x, Y )]_{x=X}]$.

\noindent {\bf Definition 2.4} ($G$-normal distribution) A
d-dimensional random vector $X = (X_1, \cdot\cdot\cdot,X_d)$ in a
sublinear expectation space $(\Omega, {\cal H}, \hat{E})$ is called
$G$-normal distributed if for each $a, b\in R$ we have $$aX +
b\hat{X}\sim \sqrt{a^2 + b^2}X,$$  where $\hat{X}$ is an independent
copy of $X$. Here the letter $G$ denotes the function $$G(A) :=
\frac{1 }{2}\hat{ E}[(AX,X)] : S_d \rightarrow R,$$  where $S_d$
denotes the collection of $d\times d$ symmetric matrices.

The function $G(\cdot) : S_d \rightarrow R$ is a monotonic,
sublinear mapping on $S_d$ and $G(A) = \frac{1 }{2}\hat{
E}[(AX,X)]\leq \frac{1 }{2}|A|\hat{ E}[|X|^2]=:\frac{1
}{2}|A|\bar{\sigma}^2$ implies that there exists a bounded, convex
and closed subset $\Gamma\subset S_d^+$ such that
$$G(A)=\frac{1 }{2}\sup_{\gamma\in \Gamma}Tr(\gamma A).$$ If there exists some $\beta>0$ such that
$G(A)-G(B)\geq \beta Tr(A-B)$ for any $A\geq B$, we call the
$G$-normal distribution is non-degenerate, which is the case we
consider throughout this article.

\noindent {\bf Definition 2.5} i) Let $\Omega_T=C_0([0, T]; R^d)$
with the supremum norm, $ {\cal H}^0_T:=\{\varphi(B_{t_1},...,
B_{t_n})| \forall n\geq1, t_1, ..., t_n \in [0, T], \forall \varphi
\in C_{l, Lip}(R^{d\times n})\}$, $G$-expectation is a sublinear
expectation defined by
$$\hat{E}[\varphi( B_{t_1}-B_{t_0}, B_{t_2}-B_{t_1} ,
\cdot\cdot\cdot, B_{t_m}- B_{t_{m-1}} )]$$$$ = \tilde{E}
[\varphi(\sqrt{t_1-t_0}\xi_1, \cdot\cdot\cdot, \sqrt{t_m
-t_{m-1}}\xi_m)],$$ for all $X=\varphi( B_{t_1}-B_{t_0},
B_{t_2}-B_{t_1} , \cdot\cdot\cdot, B_{t_m}- B_{t_{m-1}} )$, where
$\xi_1, \cdot\cdot\cdot, \xi_n$ are identically distributed
$d$-dimensional $G$-normal distributed random vectors in a sublinear
expectation space $(\tilde{\Omega}, \tilde{\cal H},\tilde{ E})$ such
that  $\xi_{i+1}$ is independent to $(\xi_1, \cdot\cdot\cdot,
\xi_i)$ for each $i = 1, \cdot\cdot\cdot,m$. $(\Omega_T, {\cal
H}^0_T, \hat{E})$ is called a $G$-expectation space.

ii) For $t\in [0, T]$ and $\xi=\varphi(B_{t_1},..., B_{t_n})\in
{\cal H}^0_T$, the conditional expectation defined by(there is no
loss of generality, we assume $t=t_i$) $$\hat{E}_{t_i}[\varphi(
B_{t_1}-B_{t_0}, B_{t_2}-B_{t_1} , \cdot\cdot\cdot, B_{t_m}-
B_{t_{m-1}} )]$$$$=\tilde{\varphi}( B_{t_1}-B_{t_0}, B_{t_2}-B_{t_1}
, \cdot\cdot\cdot, B_{t_i}- B_{t_{i-1}} ),$$ where
$$\tilde{\varphi}(x_1, \cdot\cdot\cdot, x_i)=\hat{E}[\varphi( x_1,
\cdot\cdot\cdot,x_i, B_{t_{i+1}}- B_{t_{i}}, \cdot\cdot\cdot,
B_{t_m}- B_{t_{m-1}} )].$$

Let $\|\xi\|_{p, G}=[\hat{E}(|\xi|^p)]^{1/p}$ for $\xi\in{\cal
H}^0_T$ and $p\geq1$, then $\forall t\in[0, T]$, $\hat{E}_t(\cdot)$
is a continuous mapping on ${\cal H}^0_T$ with norm $\|\cdot\|_{1,
G}$ and therefore can be extended continuously to the completion
$L^1_G(\Omega_T)$ of ${\cal H}^0_T$ under norm $\|\cdot\|_{1, G}$.

\noindent {\bf Theorem 2.6}([DHP08]) There exists a tight subset
${\cal P}\subset {\cal M}_1(\Omega_T)$ such that
$$\hat{E}(\xi)=\max_{P\in {\cal P}}E_P(\xi) \ \ \textrm{for \
all} \ \xi\in{\cal H}^0_T.$$ ${\cal P}$ is called a set that
represents $\hat{E}$.

\noindent {\bf Remark 2.7} i) Let ${\cal A}$ denotes the sets that
represent $\hat{E}$. ${\cal P}^*=\{P\in {\cal M}_1(\Omega_T)|
E_P(\xi)\leq \hat{E}(\xi), \ \forall \ \xi\in{\cal H}^0_T\}$ is
obviously the maximal one, which is  convex and weak compact. All
capacities induced by weak compact sets of probabilities in ${\cal
A}$ are the same, i.e. $c_{\cal P}:=\sup_{P\in{\cal
P}}P=\sup_{P\in{\cal P}'}P=:c_{{\cal P}'}$ for any weak compact set
${\cal P}, {\cal P}'\in{\cal A}$.

ii) Let $(\Omega^0, \{{\cal F}^0_t\}, {\cal F}, P^0 )$ be a filtered
probability space, and $\{W_t\}$ be a d-dimensional Brownian motion
under $P^0$.  [DHP08] proved that $${\cal P}'_M:=\{P_0\circ X^{-1}|
X_t=\int_0^th_sdW_s, h\in L^2_{\cal F}([0,T]; \Gamma^{1/2}) \}\in
{\cal A},$$ where $\Gamma^{1/2}:=\{\gamma^{1/2}| \gamma\in \Gamma\}$
and $\Gamma$ is the set in the representation of $G(\cdot)$.

iii) Let ${\cal P}_M$ be the weak closure of ${\cal P}'_M$. Then
under each $P\in {\cal P}_M$, the canonical process
$B_t(\omega)=\omega_t$ for $\omega\in\Omega_T$ is a martingale.

 \noindent {\bf Definition 2.8} i) Let $c$ be the capacity induced by
$\hat{E}$. A map $X$ on $\Omega_T$ with values in a topological
space is said to be quasi-continuous w.r.t $c$ if
$$\forall \varepsilon>0, \ there \ exists \ an \ open \ set \ O \ with \ c(O)<\varepsilon \ such \ that \ X|_{O^c} \ is \ continuous.$$

ii) We say that $X: \Omega_T\rightarrow R$ has a quasi-continuous
version if there exists a quasi-continuous function $Y:
\Omega_T\rightarrow R$ with $X=Y$, $c$-q.s.. $\Box$

 Let $\|\varphi\|_{p,G}=[\hat{E}(|\varphi|^p)]^{1/p}$ for $\varphi\in C_b(\Omega_T)$, the
 completions of
$C_b(\Omega_T)$, ${\cal H}^0_T$ and $L_{ip}(\Omega_T)$ under
$\|\cdot\|_{p,G}$ are the same and denoted by $L^p_G(\Omega_T)$,
where $$L_{ip}(\Omega_T):=\{\varphi(B_{t_1},..., B_{t_n})| \forall
n\geq1, t_1, ..., t_n \in [0, T], \forall \varphi \in C_{b,
Lip}(R^{d\times n})\}$$ and $C_{b, Lip}(R^{d\times n})$ denotes the
set of bounded Lipschitz functions on $R^{d\times n}$.

\noindent {\bf Theorem 2.9}[DHP08] For  $p\geq1$ the completion
$L^p_G(\Omega_T)$ of $C_b(\Omega_T)$ is $$L^p_G(\Omega_T)=\{X\in
L^0: X \ has \ a \ q.c. \ version, \
\lim_{n\rightarrow\infty}\hat{E}[|X|^p1_{\{|X|>n\}}]=0\},$$ where
$L^0$ denotes the space of all R-valued measurable functions on
$\Omega_T$.

\subsection{Basic notions on stochastic calculus
in sublinear expectation space}

\textit{For convenience of description, we only give the definition
of It$\hat{o}$ integral with respect to 1-dimensional $G$-Brownian
motion. However, all results in the following sections of this
article hold for the $d$-dimensional case.}

Let $H^0_G(0, T)$ be the collection of processes in the following
form: for a given partition $\{t_0, \cdot\cdot\cdot, t_N\} = \pi_T$
of $[0, T]$,  $$ \eta_t(\omega) = \sum^{N-1}_{j=0}
\xi_j(\omega)1_{[t_j ,t_{j+1})}(t),$$ where $\xi_i\in L_{ip}(\Omega_
{t_i})$, $i = 0, 1, 2, \cdot\cdot\cdot, N-1$. For each $\eta\in
H^0_G(0, T)$, let
$\|\eta\|_{H^{p}_G}=\{\hat{E}(\int_0^T|\eta_s|^2ds)^{p/2}\}^{1/p}$
and denote $H^{p}_G(0, T)$ the completion of $H^0_G(0, T)$ under
norm $\|\cdot\|_{H^{p}_G}$.

\noindent {\bf Definition 2.10} For each $\eta\in H^0_G(0, T)$ with
the form $$\eta_t(\omega) = \sum^{N-1}_{j=0}
\xi_j(\omega)1_{[t_j,t_{j+1})}(t),$$ we define $$I(\eta) =\int_0^T
\eta(s)dB_s := \sum^{N-1}_{j=0} \xi_j(B_{t_{j+1} }-B_{t_j} ).$$

By B-D-G inequality, the mapping $I: H^0_G(0, T)\rightarrow
L^p_G(\Omega_T)$ is continuous under $\|\cdot\|_{H^{p}_G}$ and thus
can be continuously extended to $H^p_G(0, T)$.

\noindent {\bf Definition 2.11} A process $\{M_t\}$ with values in
$L^1_G(\Omega_T)$ is called a $G$-martingale if $\hat{E}_s(M_t)=M_s$
for any $s\leq t$. If $\{M_t\}$ and  $\{-M_t\}$ are both
$G$-martingale, we call $\{M_t\}$ symmetric $G$-martingale.

\noindent {\bf Definition 2.12}  For two process $\{X_t\}, \{Y_t\}$
with values in $L^1_G(\Omega_T)$, we say $\{X_t\}$ is a version of $
\{Y_t\}$ if
$$X_t=Y_t, \ \ q.s. \ \ \forall t\in[0, T].$$

\section{ Some lemmas}

\noindent {\bf Definition 3.1}  We say that a process $\{M_t\}$ with
values in $L^1_G(\Omega_T)$ is quasi-continuous if

\emph{$\forall \varepsilon>0$, there exists open set $G$ with
$c(G)<\varepsilon$ such that $M_\cdot(\cdot)$ is continuous on
$G^c\times[0, T]$}.

\noindent {\bf Lemma 3.2}(Song10) Any $G$-martingale $\{M_t\}$  has
a quasi-continuous version. $\Box$

So we shall only consider quasi-continuous $G$-martingale in the
rest of the article. The following lemma is the counterpart of
Doob's uniform integrability lemma, and the proof is adapted from
[Yan98].

Let ${\cal B}_t=\sigma\{B_s| s\leq t\}$, ${\cal F}_t=\cap_{r>t}
{\cal B}_r$ and $\mathbb{F}=\{{\cal F}_t\}_{t\in[0, T]}$. $\tau:
\Omega_T\rightarrow[0, T]$ is called a $\mathbb{F}$ stopping time if
$[\tau\leq t]\in {\cal F}_t$, $\forall t\in[0, T]$.

\noindent {\bf Lemma 3.3} Let $\{M_t\}$ be a symmetric or negative
 $G$-martingale with $M_T\in L^p_G(\Omega_T)$ for $p\geq1$, then
$\{|M_{\sigma_i}|^p\}_{i\in I}$ are uniformly integrable under
$\hat{E}$ in the following sense:
$$\sup_{i\in
I}\hat{E}[|M_{\sigma_i}|^p1_{[|M_{\sigma_i}|>n]}]\rightarrow0,$$
where $\{\sigma_i| \ i\in I\}$ is a family of stopping times w.r.t
$\mathbb{F}$.

{\bf Proof.} Fix $P\in {\cal P}_M$ and $i\in I$.
\begin {eqnarray*}& &E_P[|M_{\sigma_i}|^p1_{[|M_{\sigma_i}|>n]}]\\
&\leq&E_P[|M_T|^p1_{[|M_{\sigma_i}|>n]}]\\
&\leq&\delta^pP(|M_{\sigma_i}|>n)+E_P[|M_T|^p1_{[|M_T|>\delta]}]\\
&\leq&\delta^pn^{-p}E_P(|M_{\sigma_i}|^p)+E_P[|M_T|^p1_{[|M_T|>\delta]}]\\
&\leq&\delta^pn^{-p}E_P(|M_T|^p)+E_P[|M_T|^p1_{[|M_T|>\delta]}].
\end {eqnarray*}

So $\sup_{i\in I}\hat{E}[|M_{\sigma_i}|^p1_{[|M_{\sigma_i}|>n]}]
\leq\delta^pn^{-p}\hat{E}(|M_T|^p)+\hat{E}[|M_T|^p1_{[|M_T|>\delta]}]$.
First let $n\rightarrow\infty$, then let $\delta$ go to infinity, we
get the result.$\Box$

\noindent {\bf Lemma 3.4} Let $E$ be a  metric space and a mapping
$E\times[0, T]\ni(\omega, t)\rightarrow M_t(\omega)\in R$ be
continuous on $E\times[0, T] $.

Define $\underline{\tau}_a=\inf\{t\geq0| \ M_t\geq a\}\wedge T$ and
$\overline{\tau}_a=\inf\{t\geq0| \ M_t> a\}\wedge T$. Then

i) $M_{t\wedge\underline{\tau}_a}$ is continuous at any $\omega\in
E$ with $M_{t\wedge\underline{\tau}_a}(\omega)<a$ and $M_{t\wedge
\overline{\tau}_a}$ is continuous at any $\omega\in E$ with
$M_{t\wedge\overline{\tau}_a}(\omega)=a$. Moreover,
$-M_{t\wedge\underline{\tau}_a}$, $M_{t\wedge\overline{\tau}_a}$ are
both lower semi-continuous.

ii)$-\overline{\tau}_a$ and $\underline{\tau}_a$ are both lower
semi-continuous.

{\bf Proof.} i) For $\omega$ with $M_{t\wedge
\overline{\tau}_a}(\omega)=a$, $M_{t\wedge
\overline{\tau}_a}(\cdot)$ is obviously continuous at $\omega$.
Also, we claim that for $\omega$ with $M_{t\wedge
\underline{\tau}_a}(\omega)<a$, $M_{t\wedge
\underline{\tau}_a}(\cdot)$ is continuous at $\omega$. Otherwise,
there exists a sequence $\{\omega_n\}\subset \Omega_T$ and a
sequence $\{t_n\}\subset[0, t]$ such that
$\omega_n\rightarrow\omega$ and $M_{t_n}(\omega_n)\geq a$. Assume
$t_n\rightarrow t'\in[0, t]$, then
\begin {eqnarray*}|M_{t_n}(\omega_n)-M_{t'}(\omega)|
\rightarrow0.
\end {eqnarray*} So $M_{t'}(\omega)\geq a$ and $M_{t\wedge
\underline{\tau}_a}(\omega)\geq a$, which contradicts the
assumption.

For any $b\in R$, we claim that $[M_{t\wedge\underline{\tau}_a}<b]$
and $[M_{t\wedge\overline{\tau}_a}>b]$ are both open. If $b>a$,
$[M_{t\wedge\underline{\tau}_a}<b]$ is obvious open. Assume $b\leq
a$. For any $\omega\in [M_{t\wedge\underline{\tau}_a}<b]$, there
exists an open set $O$ such that $\omega\in
O\subset[M_{t\wedge\underline{\tau}_a}<b]$ since
$M_{t\wedge\underline{\tau}_a}$ is continuous at $\omega$. So
$[M_{t\wedge\underline{\tau}_a}<b]$ is open. Also,
$[M_{t\wedge\overline{\tau}_a}>b]$ is obvious open for $b\geq a$.
Assume $b<a$. If $M_{t\wedge\overline{\tau}_a}(\omega)=a$, there
exists an open set $O$ such that $\omega\in
O\subset[M_{t\wedge\overline{\tau}_a}>b]$ since
$M_{t\wedge\overline{\tau}_a}$ is continuous at $\omega$. For
$b<M_{t\wedge\overline{\tau}_a}(\omega)<a$, we have $b<M_t(\omega)$.
Then there exists an open set $O$ such that $\omega\in
O\subset[M_t>b]\subset[M_{t\wedge\overline{\tau}_a}>b]$ since
$M_{t}$ is continuous at $\omega$. So
$[M_{t\wedge\overline{\tau}_a}>b]$ is open.

ii) For any $t\in[0, T]$, $[\overline{\tau}_a<t]$ is obviously open.
For any $t\in[0, T)$,
$[\underline{\tau}_a>t]=[M_{t\wedge\underline{\tau}_a}<a]$ is open
by i). $\Box$

\noindent {\bf Lemma 3.5} For any closed set $F$, we have
$$c(F)=\inf \{c(O)| \ F\subset O\},$$ where $c$ is the capacity
induced by $\hat{E}$.

{\bf Proof.} It suffices to prove that for any closed set
$F\subset\Omega_T$, $c(F)\geq\inf \{c(O)| \ F\subset O\}$. In fact,
for any closed set $F\subset\Omega_T$, there exists
$\{\varphi_n\}\in C_b(\Omega_T)$ such that $1\geq\varphi_n\downarrow
1_F$. By Theorem 28 in [DHP08], we have
$c(F)=\lim_{n\rightarrow\infty}\hat{E}(\varphi_n)$. Let
$O_n=[\varphi_n>1-1/n]$. Then $O_n\supset F$ and
$c(O_n)\leq\frac{n}{n-1}\hat{E}(\varphi_n)\rightarrow c(F)$. So
$c(F)\geq\inf_nc(O_n)\geq\inf \{c(O)| \ F\subset O\}$. $\Box$

\section{ Hitting times for $G$-martingale}

\subsection{Hitting times for symmetric $G$-martingale}

In this section, we try to define stopped processes for symmetric
$G$-martingale.

Let $${\cal Q}_T=\{(r, s)| \ T\geq r>s\geq0, \ r,s \textmd{\ are \
rational}\}$$ and $${\cal S}_a(M)=\{\omega\in \Omega_T| \ \exists
(r, s)\in {\cal Q}_T \textmd{\ such \ that }\ M_t(\omega)=a \
\forall t\in [s, r]\}.$$

\noindent {\bf Theorem 4.1} Let $\{M_t\}_{t\in[0, T]}$ be a
 symmetric $G$-martingale. Then for all $a>M_0$  and
$\overline{\tau}_a, \underline{\tau}_a$ defined above,

i)$\forall t\in[0, T]$, $M_{t\wedge \overline{\tau}_a}$ and
$M_{t\wedge \underline{\tau}_a}$ are both quasi-continuous.
Consequently, $\{M_{t\wedge \overline{\tau}_a}\}$ and $\{M_{t\wedge
\underline{\tau}_a}\}$ are both symmetric $G$-martingale.

ii) If in addition $c({\cal S}_a(M))=0$, then $\overline{\tau}_a,
\underline{\tau}_a$ are both quasi-continuous.

{\bf Proof.} i) Since $\{M_t\}_{t\in[0,T]}$ be a symmetric
$G$-martingale, it is a martingale under each $P\in {\cal P}_M$.
Therefore, $E_P(M_{t\wedge \overline{\tau}_a})=M_0=E_P(M_{t\wedge
\underline{\tau}_a})$ for each $P\in {\cal P}_M$. Consequently
$\hat{E}(M_{t\wedge \underline{\tau}_a}-M_{t\wedge
\overline{\tau}_a})=0$. Noting that $M_{t\wedge
\underline{\tau}_a}\geq M_{t\wedge \overline{\tau}_a}$, we get
$M_{t\wedge \underline{\tau}_a}= M_{t\wedge \overline{\tau}_a}$,
q.s.  Since $\{M_t\}$ is quasi-continuous, for any $\varepsilon>0$
there there exists open set $G$ with $c(G)<\varepsilon/2$ such that
$M_\cdot(\cdot)$ is continuous on $G^c\times[0, T]$. Let ${\cal
Q}=\{(r, s)| \ r>s, \ r,s \textmd{\ are \ rational}\}$. Noting that
$$[M_{t\wedge \underline{\tau}_a}> M_{t\wedge \overline{\tau}_a}]=\cup_{(r,s)\in{\cal
Q} }[M_{t\wedge \underline{\tau}_a}\geq r, s\geq M_{t\wedge
\overline{\tau}_a}],$$ we have
$$[M_{t\wedge \underline{\tau}_a}> M_{t\wedge \overline{\tau}_a}]\subset G\bigcup\cup_{(r,s)\in{\cal
Q} }([M_{t\wedge \underline{\tau}_a}\geq r, s\geq M_{t\wedge
\overline{\tau}_a}]\cap G^c).$$ By Lemma 3.4, $[M_{t\wedge
\underline{\tau}_a}\geq r, s\geq M_{t\wedge \overline{\tau}_a}]\cap
G^c$ is closed for any $(r, s)\in {\cal Q}$. Since $c([M_{t\wedge
\underline{\tau}_a}\geq r, s\geq M_{t\wedge \overline{\tau}_a}]\cap
G^c)=0$, by Lemma 3.5 there exists open set $O$ with
$c(O)<\varepsilon/2$ such that $$\cup_{(r,s)\in{\cal Q}
}([M_{t\wedge \underline{\tau}_a}\geq r, s\geq M_{t\wedge
\overline{\tau}_a}]\cap G^c)\subset O.$$  By Lemma 3.3, $M_{t\wedge
\overline{\tau}_a}$ and $M_{t\wedge \underline{\tau}_a}$ are both
continuous on $O^c\cap G^c$.

ii) By the quasi-continuity of $\{M_t\}$, for any $\varepsilon>0$
there exists open set $G$ such that $c(G)<\varepsilon/2$ and
$M_t(\omega)$ is continuous on $G^c\times[0, T]$. So
$$G^c\cap[\overline{\tau}_a>\underline{\tau}_a]\subset{\cal S}_a(M)\bigcup\cup_{r\in
\mathbb{Q}\cap[0, T]}[M_{r\wedge \overline{\tau}_a}<M_{r\wedge
\underline{\tau}_a}],$$ where $\mathbb{Q}$ denotes the totality of
rational numbers. Then
$c(G^c\cap[\overline{\tau}_a>\underline{\tau}_a])=0$. Since
$$G^c\cap[\overline{\tau}_a>\underline{\tau}_a]=\bigcup_{(r,s)\in{\cal Q}
}([\overline{\tau}_a\geq r, s\geq \underline{\tau}_a]\cap G^c)$$ and
$[\overline{\tau}_a\geq r, s\geq \underline{\tau}_a]\cap G^c$ is
closed by Lemma 3.4, there exists open set $O$ such that
$c(O)<\varepsilon/2$ and
$G^c\cap[\overline{\tau}_a>\underline{\tau}_a]\subset O$. So on
$O^c\cap G^c$, $\overline{\tau}_a=\underline{\tau}_a$ are both
continuous. $\Box$

\noindent {\bf Remark 4.2} If the quadratic variation process of
$\{M_t\}$ is strictly increasing except on a polar set, then
$c({\cal S}_a(M))=0$ for any $a\in \mathbb{R}$.

\noindent {\bf Example 4.3} Let $\{B_t\}_{t\in[0, T]}$ be a
1-dimensional  $G$-Brownian motion. For $a>0$, let
$\underline{\tau}_a=\inf\{t\geq0| \ B_t\geq a\}\wedge T$ and
$\overline{\tau}_a=\inf\{t\geq0| \ B_t> a\}\wedge T$. Then we have

i) For any $t\in[0, T]$, $-B_{t\wedge\underline{\tau}_a}$,
$B_{t\wedge\overline{\tau}_a}$, $-\overline{\tau}_a$ and
$\underline{\tau}_a$ are all lower semi-continuous.

ii) For any $t\in[0, T]$, $B_{t\wedge\underline{\tau}_a}$,
$B_{t\wedge\overline{\tau}_a}$, $\overline{\tau}_a$ and
$\underline{\tau}_a$ are all quasi-continuous.

iii) $\{B_{t\wedge\underline{\tau}_a\}}$ and
$\{B_{t\wedge\overline{\tau}_a\}}$ are both symmetric
$G$-martingale.

\subsection{Hitting times for  $G$-martingale(non-symmetric)}

For each $P\in{\cal P}_M$ and $t\in[0, T]$, let ${\cal A}_{t,
P}:=\{Q\in{\cal P}_M| \ Q=P|_{{\cal F}_t}\}$. Theorem 2.3 in [STZ09]
implies the following result: For $t\in[0, T]$ and $\xi\in
L^1_G(\Omega_T)$, $\eta\in L^1_G(\Omega_t)$, $\eta=\hat{E}_t(\xi)$
if and only if for each $P\in{\cal P}_M$
$$\eta={\textmd{ess}\sup}_{Q\in {\cal A}_{t,
P}}^PE_Q(\xi|{\cal F}_t), \ P-a.s.$$

\noindent {\bf Theorem 4.4} Let $\{M_t\}_{t\in[0, T]}$ be a
quasi-continuous $G$-martingale.  For all $a>|M_0|$,  $M_{t\wedge
\overline{\sigma}_a}$ and $M_{t\wedge \underline{\sigma}_a}$ are
both $G$-martingale, where $\underline{\sigma}_a=\inf\{t\geq0| \
M_t\leq -a\}\wedge T$ and $\overline{\sigma}_a=\inf\{t\geq0| \ M_t<
-a\}\wedge T$.

{\bf Proof.} For each $P\in{{\cal P}_M}$, $\{M_t\}_{t\in[0, T]}$ is
a supermartingale, so for each $t\in[0, T]$
$$E_P(M_{t\wedge\overline{\sigma}_a}|{\cal
F}_{t\wedge\underline{\sigma}_a})\leq
M_{t\wedge\underline{\sigma}_a}$$ by Doob optimal stopping theorem
and  noting that $\underline{\sigma}_a\leq \overline{\sigma}_a$.
This implies $E_P(M_{t\wedge\underline{\sigma}_a}-
M_{t\wedge\overline{\sigma}_a})\geq0.$ On the other hand, it's
obvious to see that $M_{t\wedge\underline{\sigma}_a}\leq
M_{t\wedge\overline{\sigma}_a}$. So
$M_{t\wedge\underline{\sigma}_a}= M_{t\wedge\overline{\sigma}_a}$
q.s. By the same arguments as in Theorem 4.1, for any
$\varepsilon>0$, there exists an open set $O$ such that
$c(O)<\varepsilon$ and $M_{t\wedge\underline{\sigma}_a}=
M_{t\wedge\overline{\sigma}_a}$ are continuous on $O^c$.  So
$M_{t\wedge\underline{\sigma}_a}$ and
$M_{t\wedge\overline{\sigma}_a}$ are both quasi continuous.

Let $\sigma=\overline{\sigma}_a$ or $\underline{\sigma}_a$.

For $0\leq s<t\leq T$ and $P\in{\cal P}_M$,
\begin {eqnarray*}& &\hat{E}_s(M_{t\wedge\sigma})\\
&=&{\textmd{ess}\sup}_{Q\in {\cal A}_{s,
P}}^PE_Q(M_{t\wedge\sigma}|{\cal F}_s)\\
&\leq& M_{s\wedge\sigma} \ \ \ \ P-a.s.
\end {eqnarray*}

On the other hand,
\begin {eqnarray*}
& &E_Q(M_{t\wedge\sigma}|{\cal F}_s)\\
&=&-a1_{[\sigma\leq s]}+E_Q(M_{t\wedge\sigma}|{\cal F}_s)1_{[\sigma>
s]}\\
&\geq&-a1_{[\sigma\leq s]}+E_Q[E_Q(M_t|{\cal
F}_{t\wedge\sigma})|{\cal F}_s]1_{[\sigma> s]}\\
&=& -a1_{[\sigma\leq s]}+E_Q(M_t|{\cal
F}_{s\wedge\sigma})1_{[\sigma> s]}\\
&=&-a1_{[\sigma\leq s]}+E_Q(M_t|{\cal F}_s)1_{[\sigma> s]}.
\end {eqnarray*}

So
\begin {eqnarray*}& &\hat{E}_s(M_{t\wedge\sigma})\\
&=&{\textmd{ess}\sup}_{Q\in {\cal A}_{s,
P}}^PE_Q(M_{t\wedge\sigma}|{\cal F}_s)\\
&\geq& -a1_{[\sigma\leq s]}+{\textmd{ess}\sup}_{Q\in {\cal A}_{s,
P}}^PE_Q(M_t|{\cal F}_s)1_{[\sigma> s]}\\
&=& -a1_{[\sigma\leq s]}+M_s1_{[\sigma> s]}\\
&=& M_{\sigma\wedge s} \ \ P-a.s.
\end {eqnarray*}
$\hat{E}_s(M^\sigma_t)= M^\sigma_s$ q.s. $\Box$

\section {Applications}

\noindent {\bf Theorem 5.1} Let $\xi\in L^\beta(\Omega_T)$ for some
$\beta\geq1$ be symmetric, then there exist a sequence
$\{\xi^n\}\subset L^1(\Omega_T)$ which are bounded and symmetric
such that $\hat{E}[|\xi-\xi^n|^\beta]\rightarrow0$.

{\bf Proof.} Let $M_t=\hat{E}_t(\xi)$ for $t\in[0, T]$ be the
quasi-continuous version. For each $n\in N$, let
$\sigma_n=\inf\{t\geq0| \ |M_t|> n\}\wedge T$, and
$\tau_n=\inf\{t\geq0| \ M_t> n\}\wedge T$. By Theorem 4.1,
$\{M^{\tau_n}_t\}_{t\in[0, T]}$ is a symmetric $G$-martingale. Let
$\{N_t\}$ be the quasi-continuous version of $M^{\tau_n}$  and
$\varsigma_n=\inf\{t\geq0| \ -N_t>n\}\wedge T$. By the same
arguments, $\{N_t^{\varsigma_n}\}$ is a bounded symmetric
$G$-martingale. Since the paths of $\{M^{\tau_n}_t\}$ and
$\{N_t^{\varsigma_n}\}$ are continuous except on a polar set,
$$\{\omega\in\Omega_T| \exists \textmd{\ some} \ t\in[0, T], \ s.t. \ N_t(\omega)\neq M^{\tau_n}_t(\omega)
\}$$ is a polar set. So
$N_t^{\varsigma_n}=M^{\tau_n\wedge\varsigma_n}_t=M^{\sigma_n}_t, \
q.s..$
\begin {eqnarray*}
& &|M_{\sigma_n}-M_T|^\beta\\
&\leq& 2^{\beta-1}(|M_{\sigma_n}-(M_T\wedge n)\vee(-n)|^\beta+|(M_T\wedge n)\vee(-n)-M_T|^\beta)\\
&\leq&2^{2\beta-1}|M_{\sigma_n}|^\beta1_{[|M_{\sigma_n}|\geq
n]}+2^{\beta-1}|M_T|^\beta1_{[|M_T|>n]}.
\end {eqnarray*}
Hence, by Lemma 3.3,
$$\hat{E}(|M_{\sigma_n}-M_T|^\beta)\leq 2^{2\beta-1}\sup_i\hat{E}[|M_{\sigma_i}|^\beta1_{[|M_{\sigma_i}|\geq n]}]+2^{\beta-1}\hat{E}[|M_T|^\beta1_{[|M_T|>n]}]\rightarrow0.$$
$\Box$

The following Corollary improved Theorem 4.6 in [Song10].

\noindent {\bf Corollary 5.2} Let $\xi\in L^\beta_G(\Omega_T)$ for
some $\beta>1$ with $\hat{E}(\xi)+\hat{E}(-\xi)=0$, then there
exists $\{Z_t\}_{t\in[0, T]}\in H^\beta_G(0, T)$ such that
$$\xi=\hat{E}(\xi)+\int_0^TZ_sdB_s.$$

{\bf Proof.} By Theorem 5.1, there exist a sequence
$\{\xi^n\}\subset L^1(\Omega_T)$ which are bounded and symmetric
such that $\hat{E}[|\xi-\xi^n|^\beta]\rightarrow0$. By Theorem 4.6
in [Song10], there exists $\{Z^n_t\}_{t\in[0, T]}\in H^\beta_G(0,
T)$ such that
$$\xi^n=\hat{E}(\xi^n)+\int_0^TZ^n_sdB_s.$$ By B-D-G and Doob's maximal
inequality, $\{Z^n_t\}_{t\in[0, T]}$ is a Cauchy  sequence in
$H^\beta_G(0, T)$.  So there exists $\{Z_t\}\in H^\beta_G(0, T)$
such that $\|Z^n-Z\|_{H^\beta_G}\rightarrow0$. Then
$$\xi=\lim_{L^\beta_G, n\rightarrow\infty}\xi^n=\lim_{L^\beta_G, n\rightarrow\infty}[\hat{E}(\xi^n)+\int_0^TZ^n_sdB_s]=\hat{E}(\xi)+\int_0^TZ_sdB_s.$$
$\Box$

\noindent {\bf Theorem 5.3} Let $\xi\in L^2_G(\Omega_T)$ be bounded
above, then $M_t=\hat{E}_t(\xi)$, $ t\in[0, T]$ has the following
representation:
\begin {eqnarray}
M_t=M_0+\int_0^tZ_sdB_s-K_t,
\end {eqnarray}
 where $\{Z_t\}\in M^2_G(\Omega_T)$, $\{K_t\}$ is a quasi-continuous
 increasing process with $K_0=0$ and $\{-K_t\}_{t\in[0, T]}$ a
$G$-martingale.

{\bf Proof.} There is no loss of generality, we only consider the
$\xi\leq0$ case. Let $\{M_t\}$ be the quasi-continuous version. For
$n\in N$, let $\sigma_n=\overline{\sigma}_n$ defined in Theorem 4.4.
Then $M_{\sigma_n}$ is bounded. By Theorem 4.4 and Theorem 4.5 in
[Song10], we have the following representation
\begin {eqnarray*}
M^{\sigma_n}_t=\hat{E}(M_T)+\int_0^tZ^n_sdB_s-K^n_t=:N^n_t-K^n_t, \
q.s.
\end {eqnarray*}
where $\{Z^n_t\}_{t\in[0, T]}\in H^2_G(0, T)$ and
$\{K^n_t\}_{t\in[0, T]}$ is a continuous increasing process
 with $K^n_0=0$ and $\{-K^n_t\}_{t\in[0, T]}$ a
$G$-martingale. For $m>n$, by uniqueness of decomposition of
semimartingale, $N^n_t=(N^m)^{\sigma_n}_t$ and
$K^n_t=(K^m)^{\sigma_n}_t$. So
$$\widehat{M}_t:=M^{\sigma_m}_t-M^{\sigma_n}_t=(N^m_t-N^n_t)-(K^m_t-K^n_t)=:\widehat{N}_t-\widehat{K}_t$$
with $\{\widehat{K}_t\}$ a continuous increasing process.

By It$\hat{o}$ formula, we have
$$\hat{E}(\widehat{N}^2_T)\leq\hat{E}(\widehat{M}^2_T)+2\hat{E}(\int_0^T\widehat{M}^+_sd\widehat{K}_s).$$
Noting that $\widehat{M}^+_s\leq n$, we have
$$\hat{E}(\widehat{N}^2_T)\leq\hat{E}(\widehat{M}^2_T)+2n\hat{E}(\widehat{K}_T).$$

$\hat{E}(\widehat{M}^2_T)\leq2\{\hat{E}[(M_T-M_{\sigma_n})^2]+\hat{E}[(M_T-M_{\sigma_m})^2]\}$,
\begin {eqnarray*}& &2n\hat{E}(\widehat{K}_T)\\
&=&2n[\hat{E}(M_{\sigma_n}-M_{\sigma_m})+\hat{E}(M_{\sigma_m}-M_{\sigma_n})]\\
&\leq&4n[\hat{E}(|M_{\sigma_n}-M_T|)+\hat{E}(|M_{\sigma_m}-M_T|)].
\end {eqnarray*}
By the same arguments as in Theorem 5.1,
$\hat{E}(\widehat{M}^2_T)\rightarrow0$ as $m, n\rightarrow\infty$.

\begin {eqnarray*}
& &|M_{\sigma_n}-M_T|\\
&\leq& |M_{\sigma_n}-M_T\vee(-n)|+|M_T\vee(-n)-M_T|\\
&\leq&2M_{\sigma_n}1_{[|M_{\sigma_n}|\geq n]}+M_T1_{[|M_T|>n]}.
\end {eqnarray*}
So
\begin {eqnarray*}& &2n\hat{E}(\widehat{K}_T)\\
&\leq&8n\hat{E}(M_{\sigma_n}1_{[|M_{\sigma_n}|\geq
n]})+4n\hat{E}(M_T1_{[|M_T|>n]})+\\
& &8m\hat{E}(M_{\sigma_m}1_{[|M_{\sigma_m}|\geq
m]})+4m\hat{E}(M_T1_{[|M_T|>m]})\\
&\leq&8\hat{E}(M_{\sigma_n}^21_{[|M_{\sigma_n}|\geq
n]})+4\hat{E}(M_T^21_{[|M_T|>n]})+\\
& &8\hat{E}(M_{\sigma_m}^21_{[|M_{\sigma_m}|\geq
m]})+4\hat{E}(M_T^21_{[|M_T|>m]}).
\end {eqnarray*}
So $2n\hat{E}(\widehat{K}_T)\rightarrow0$ as $m,n\rightarrow\infty$
by Lemma 3.3. Consequently, we conclude that
$\hat{E}(\widehat{N}^2_T)\rightarrow0$ and
$\hat{E}(\widehat{K}^2_T)\rightarrow0$ as $m,n\rightarrow\infty$.

So $\{Z_t^n\}$ and $\{K^n_t\}$ be Cauchy sequences in $H^2_G(0, T)$
and $L^2_G(\Omega_T)$ respectively. Let $\{Z_t\}_{t\in[0, T]}$,
$\{K_t\}$ be the corresponding limits of $\{Z^n_t\}_{t\in[0, T]}$,
$\{K^n_t\}$. Then
$$M_t=\lim_{L^2_G,n\rightarrow\infty}M^n_t=\lim_{L^2_G,n\rightarrow\infty}\int_0^tZ^n_sdB_s-\lim_{L^2_G,n\rightarrow\infty}K^n_t=\int_0^tZ_sdB_s-K_t.$$
$\Box$

In this theorem, for $\xi\in L^2_G(\Omega_T)$ and bounded above, we
have $K_T\in L^2_G(\Omega_T)$. In this sense, this result improved
Theorem 4.5 in [Song10].



\providecommand{\bysame}{\leavevmode\hbox
to3em{\hrulefill}\thinspace}
\providecommand{\MR}{\relax\ifhmode\unskip\space\fi MR }
\providecommand{\MRhref}[2]{%
  \href{http://www.ams.org/mathscinet-getitem?mr=#1}{#2}
} \providecommand{\href}[2]{#2}

\end{document}